\scrollmode
\documentclass[12pt]{amsart}
\usepackage{amssymb}
\usepackage{verbatim}
\usepackage{eucal}
\usepackage{mathrsfs}
\usepackage{eepic}
\usepackage{graphicx}
\usepackage{psfrag}
\theoremstyle{plain}
\newtheorem{theorem}[subsection]{Theorem}
\newtheorem{proposition}[subsection]{Proposition}

\theoremstyle{definition}
\newtheorem{definition}[subsection]{Definition}

\theoremstyle{remark}

\newcommand{\R}{{\mathbb R}}

\newcommand{\C}{{\mathbb C}}
\newcommand{\Z}{{\mathbb Z}}
\newcommand{\Q}{{\mathbb Q}}

\newcommand{\UU}{\mathbf{U}}
\newcommand{\Proj}{{\mathbb P}}
\newcommand{\aaa}{{\mathfrak{a}}}
\newcommand{\hhh}{{\mathfrak{h}}}

\newcommand{\WWW}{{\mathscr W}}

\newcommand{\GL}{\operatorname{GL}}
\newcommand{\SL}{\operatorname{SL}}
\newcommand{\SO}{\operatorname{SO}}

\newcommand{\ww}{\mathbf{w}}
\newcommand{\cusp}{\operatorname{cusp}}
\begin{document}

\title{Eisenstein series twisted by modular symbols for the group $\SL _{n}$}
\author{Dorian Goldfeld}
\address{Department of Mathematics\\
Columbia University\\
New York, New York  10027}

\email{goldfeld@columbia.edu}

\author{Paul E. Gunnells}
\address{Department of Mathematics\\
Rutgers University\\
Newark, New Jersey 07102}

\email{gunnells@andromeda.rutgers.edu}

\thanks{Both authors were partially supported by the NSF}

\begin{abstract}
We define Eisenstein series twisted by modular symbols for the group
$\SL _{n}$, generalizing a construction of the first author
\cite{goldfeld1, goldfeld2}.  We show that, in the case of series
attached to the minimal parabolic subgroup, our series converges for
all points in a suitable cone.  We conclude with examples for $\SL
_{2}$ and $\SL _{3}$.
\end{abstract}

\date{April 18, 2000.  Revised October 1, 2000}
\maketitle

\section{Introduction}

\subsection{}
Let $\Gamma$ denote a finitely generated discrete subgroup of
$\SL_{2}(\R)$ that contains translations and acts on the upper
halfplane $\hhh$. An automorphic form of real weight $r$ and
multiplier $\psi\colon \Gamma \to \UU$ (here $\UU = \bigl\{ w \in
\C\bigm | |w| = 1\big \}$ is the unit circle) is a meromorphic
function $F\colon \hhh \to \C$ that satisfies
$$F(\gamma z) = \psi(\gamma) \cdot j(\gamma, z)^r\cdot F(z)$$
for all $\gamma =
\left(\begin{smallmatrix}a&b\\c&d\end{smallmatrix}\right)\in\Gamma$
with $j(\gamma, z) = cz+d$.  For $r \ge 0$, an integer, and $F\colon
\hhh \to \C$ any function with sufficiently many derivatives, G. Bol
\cite{bol} proved the identity
\[
\frac{d^r}{dz^{r+1} }{\Bigl ( (cz+d)^r F(\gamma z)\Bigr ) =
(cz+d)^{-r-2} F^{(r+1)} (\gamma z)},
\]
which holds for all $\gamma \in \SL_{2}(\R)$. It follows that if
$f(z)$ is an automorphic form of weight $r + 2$ and multiplier $\psi$,
and if $F$ is any $(r + 1)$-fold indefinite integral of $f$, then $F$
satisfies the functional equation
$$F(\gamma z) = \psi(\gamma) (cz+d)^{-r} \Big (F(z) \; + \;
\phi(\gamma, z)\Big ),$$
where $\phi(\gamma, z)$ is a polynomial in $z$ of degree $\le r$
satisfying the cocycle condition
$$\phi(\gamma_1 \gamma_2, z) = \overline{\psi(\gamma_2)}\, j(\gamma_2,
z)^r\, \phi(\gamma_1, \gamma_2 z) \;\; + \;\; \phi(\gamma_1, z).$$
Such a function $F$ is called an \emph{automorphic} (or
\emph{Eichler}) \emph{integral}, and the corresponding polynomial
$\phi(\gamma, z)$ is called a \emph{period polynomial}.

\subsection{}
The Eisenstein series $E^*(z; r, \psi, \phi )$ (twisted by a period
polynomial $\phi$) is defined by the infinite series
$$E^* (z; r, \psi, \phi ) = \sum_{\gamma \in \Gamma_\infty\backslash
\Gamma} \psi(\gamma) \,\phi(\gamma, z) \, j(\gamma, z)^{-r}.  $$ The
twisted Eisenstein series $E^* (z;r,\psi, \phi) $ was
first introduced by Eichler \cite{eichler} (1965).  Automorphic
integrals and Eisenstein series twisted by period polynomials were
systematically studied by Knopp \cite{knopp} (1974).  More recently
\cite{goldfeld1, goldfeld2} (1995) nonholomorphic Eisenstein series
twisted by modular symbols (period polynomials of degree 0) were
introduced (cf. \S\ref{fivethree}). O'Sullivan \cite{cormac} found
(using Selberg's method) the functional equation of these twisted
Eisenstein series; very recently O'Sullivan and Chinta
\cite{cormac.gautam} explicitly computed the scattering matrix
occurring in the functional equation.

In this paper we show how to generalize the construction of Eisenstein
series twisted by modular symbols to the group $\SL_{n}$.  The basic
properties and region of absolute convergence of such series are
obtained in the case of the minimal parabolic subgroup.  We conjecture that
these series satisfy functional equations.  

\section{Eisenstein series}\label{notation}

\subsection{}
In this section we recall the definition of cuspidal Eisenstein series
following Langlands \cite[Ch. 4]{langlands}.  We begin with some notation.  

Let $G = \SL _{n} (\R)$, let $K=\SO_{n} (\R)$, and let $\Gamma \subset
G (\Z )$ be an arithmetic group.  Let $P_{0}\subset G$ be the subgroup
of upper-triangular matrices, and let $A_{0}\subset P_{0}$ be the
subgroup of diagonal matrices with each entry positive.  For each
decomposition $n=n_{1}+\cdots+n_{k}$ with $n_{i}>0$, we have a
standard parabolic subgroup
\[
P = \Biggl\{\Biggl(\begin{array}{ccc}
P_{1}&\cdots&*\\
&\ddots&\vdots\\
0&&P_{k}
\end{array} \Biggr)\Biggm| P_{i}\in \GL _{n_{i}} (\R), \prod \det (P_{i}) = 1\Biggr\}. 
\]
We fix a Langlands decomposition $P=M_{P}A_{P}N_{P}$ as follows:
$M_{P}$ is the subgroup of block diagonal matrices, with each block
an element of $\SL ^{\pm }_{n_{i}} (\R)$; $A_{P}\subset P$ is the
subgroup with the $i$th block of the form $a_{i}I_{n_{i}}$, where
$a_{i}>0$ and $I_{n_{i}}$ is the $n_{i}\times n_{i}$ identity; and
$N_{P}\subset P$ is the subgroup with the $i$th block equal to
$I_{n_{i}}$. We transfer these decompositions to all rational
parabolic subgroups by conjugation.

\subsection{}
Let $\aaa _{0}$, $\aaa _{P}$ be the Lie algebras of the groups
$A_{0}$, $A_{P}$.  Let $\check{\aaa} _{0}$, $\check{\aaa }_{P}$ be their
$\R$-duals, and denote the pairing by $\langle
\phantom{a},\phantom{a}\rangle$.  Let $R= R^{+}\cup R^{-}\subset
\check{\aaa} _{0}$ be the roots of $G$, and let $\Delta \subset R^{+}$
be the standard set of simple roots.  For any root $\alpha $, let
$\check{\alpha }$ be the corresponding coroot.  For any parabolic
subgroup $P$, let $\rho _{P} = 1/2\sum _{\alpha \in R^{+}\cap
\check{\aaa }_{P}} \alpha $.

We recall the definition of the \emph{height function} $H_{P}\colon P\rightarrow \aaa _{P}$.  Given $p\in P$, write $p=man$, where
$m\in M_{P}$, $a\in A_{P}$, and $n\in N_{P}$.  Then $H_{P} (p)$ is defined
via
\[
e^{\langle \chi , H_{P} (p)\rangle} = a^{\chi }, \quad \text{for all
$\chi \in \check{\aaa }_{P}$.}
\]
Using an Iwasawa decomposition $G = PK$, we extend the
height function to a map $H_{P}\colon G\rightarrow \aaa _{P}$ by
setting $H_{P} (g) = H_{P} (p)$, where $g = pk$, $p\in P$, $k\in
K$.

\subsection{}\label{eisen}
Fix a parabolic subgroup $P$, and let $\Gamma _{P} = \Gamma \cap P$.
Let $f\in C^{\infty } (A_{P}N_{P}\backslash G)$ be a $\Gamma
_{P}$-invariant, $K$-finite function such that for each $g\in G$, the
function $m\mapsto f (mg)$, $m\in M_{P}$, is a square-integrable
automorphic form on $M_{P}$ with respect to $\Gamma _{P}\cap M_{P}$.
Let $\lambda \in (\check{\aaa }_{P})\otimes \C $, and let $g\in G$.
\begin{definition}
The 
\emph{Eisenstein series} associated to the above data is
\[
E_{P} (f,\lambda, g ) = \sum _{\gamma \in \Gamma _{P}\backslash \Gamma }
e^{\langle \rho _{P}+\lambda , H_{P} (\gamma g)\rangle} f (\gamma g).
\]
\end{definition}

It is known \cite[Lemma 4.1]{langlands} that this series converges
absolutely and uniformly on compact subsets of $G\times C$, where 
\begin{equation}\label{coneeq}
C = \bigl\{\lambda \bigm | \langle \Re \lambda , \check{\alpha }\rangle >
\langle \rho _{P}, \check{\alpha} \rangle, \hbox{for all $\alpha \in \Delta
$}\bigr\};
\end{equation}
here $\Re$ denotes real part. 

\section{Modular Symbols}\label{modsym}

\subsection{}
We recall the definition of modular symbols.  Our definition is
equivalent to that of Ash \cite{ash} and Ash-Borel \cite{ash-borel},
but we need a slightly different formulation for our purposes.  

Let $V = \Q ^{n}$ with the canonical $G (\Q )$-action.  Let $\ww$ be a
tuple of subspaces $(W_{1},\dots ,W_{k})$, where $W_{i}\subset V$.
The \emph{type} of $\ww$ is the tuple $(\dim W_{1},\dots ,\dim
W_{k})$.  The tuple $\ww$ is called \emph{full} if $\sum \dim W_{i} =
n$, and is called a \emph{splitting} if $V = \bigoplus _{i}
W_{i}$.  Any splitting determines a rational flag
\[
F_{\ww} = \bigl\{ \{0 \}\subsetneq F_{1}\subsetneq \cdots
\subsetneq F_{k} \subsetneq V   \bigr\}
\]
by $F_{j} = \bigoplus _{i\leq j} W_{i}$, and thus determines a
rational parabolic subgroup $P_{\ww}$, the stabilizer of $F_{\ww}$.  We abuse
notation slightly and write $P_{\ww} = M_{\ww} A_{\ww} N_{\ww}$ for
the associated Langlands decomposition.  It is easy to check that,
with our fixed decomposition, $M_{\ww} (\Q )$ preserves each $W_{i}$.

\subsection{}
Let $X$ be the symmetric space $G/K$, and let $\bar X$ be the
bordification of $X$ constructed by Borel-Serre \cite{borel-serre}.
Then the cohomology $H^{i} (\Gamma ;\C )$ may be identified with
$H^{i} (\Gamma \backslash X; \C )$ and $H^{i} (\Gamma \backslash \bar
X; \C )$.

Let $Y = \Gamma \backslash X$, $\bar Y = \Gamma \backslash \bar X$,
$\partial \bar Y = \bar Y\smallsetminus Y$, and let $\pi \colon X
\rightarrow Y$ be the canonical projection.  Let $d = (n^{2}+n)/2-1$
be the dimension of $Y$.  For all $i$, Lefschetz duality gives an
isomorphism
\[
H_{d-i} (\bar Y, \partial\bar Y;
\C)\longrightarrow H^{i} (\Gamma ; \C ).
\]

\subsection{}
Let $\ww$ be a splitting, and let $K_{\ww}$
be $K\cap M_{\ww}A_{\ww }$.  The inclusion $M_{\ww}A_{\ww }
\rightarrow G$ induces a proper map 
\[
\iota \colon M_{\ww}A_{\ww }/K_{\ww} \longrightarrow X.
\]
Let $Y_{\ww}$ be the closure of $(\pi \circ \iota) (M_{\ww}A_{\ww }/K_{\ww})$, and let $d (\ww )$ be the dimension of $Y_{\ww }$. 

\begin{definition}
Let $\ww = (W_{1},\dots ,W_{k})$ be a full tuple of subspaces.  Then
the \emph{modular symbol} $\Xi _{\ww }$ associated to $\ww $ is
defined as follows:
\begin{enumerate}
\item If $\ww $ is a splitting, then $\Xi _{\ww } \in H_{d (\ww )}
(\bar Y,
\partial \bar Y;\C)$ is the fundamental class of $Y_{\ww }$.
\item Otherwise, $\Xi _{\ww }$ is defined to be $0\in H_{d (\ww )}
(\bar Y, \partial \bar Y;\C)$, where $d (\ww )$ is the homological
degree determined by any splitting with the same type as $\ww $.
\end{enumerate}
\end{definition}

\subsection{}\label{threefive}
We define a $G (\Q )$-action on tuples as follows.  Given a full tuple
$\ww = (W_{1},\dots ,W_{k})$, let $g\cdot \ww$ be the tuple
$(W_{1},gW_{2},\dots ,gW_{k})$.  By abuse of notation we write $g\cdot
\Xi_{\ww } $ for the modular symbol $\Xi _{g\cdot \ww }$.  

Note that this is not a $G (\Q )$-action on modular symbols, since
associativity does not hold.  However, the definition $g \cdot \Xi
_{\ww }$ will suffice for our construction.

Note also that $g\cdot \Xi_{\ww }$ is different from the modular
symbol obtained via the natural $G (\Q )$-action defined by left
translation of all subspaces in a tuple.  In particular, let $\gamma
\in \Gamma $, and let $\ww ' = (\gamma W_{1},\dots ,\gamma W_{k})$.
Then $\Xi _{\ww } = \Xi _{\ww '}$, but $\Xi _{\ww } \not =\gamma \cdot
\Xi _{\ww }$ in general.

\section{Eisenstein series twisted by modular symbols}\label{}

\subsection{}
Let $\ww = (W_{1},\dots ,W_{k})$ be a full tuple of subspaces, and let
$P$ be a rational parabolic subgroup.  We say that $P$ and $\ww $ are
\emph{compatible} if the following conditions hold: there is a
splitting $\ww ' = (W_{1}',\dots ,W_{k}')$ such that $P = P_{\ww '}$,
the types of $\ww $ and $\ww '$ are equal, and $W_{1} = W_{1}'$.

Fix a rational parabolic subgroup $P$ and a compatible splitting $\ww$.  Let
$f$, $\lambda $ be as in \S\ref{eisen}, and let $\varphi $ be a $\C
$-valued linear form on $H_{d (\ww )} (\bar Y, \partial \bar Y; \C )$.  

\begin{definition}\label{tes.def}
The \emph{twisted Eisenstein series}
associated to the above data is 
\begin{equation}\label{tes}
E_{P,\varphi }^{*} = E_{P,\varphi }^{*} (f, \lambda , g, \ww) = \sum
_{\gamma \in \Gamma _{P}\backslash \Gamma } \varphi(\gamma \cdot \Xi
_{\ww})e^{\langle \rho _{P}+\lambda , H_{P} (\gamma g)\rangle} f
(\gamma g).
\end{equation}

\end{definition}

We refer to \S\ref{examples} for examples of this series, and for a
comparison with the construction in \cite{goldfeld1, goldfeld2}.

\begin{proposition}\label{wd}
The series in \eqref{tes} is well-defined.
\end{proposition}

\begin{proof}
Let $\ww = (W_{1},\dots ,W_{k})$ and let $\gamma \in \Gamma $.  We
need to show that the modular symbol $\gamma \cdot \Xi_{\ww } $
depends only on the coset $\Gamma _{P}\gamma $.

First we assume $\gamma \cdot \ww $ is a splitting.  By the remarks at
the end of \S\ref{threefive}, if $\gamma \in \Gamma $ and $\ww' =
(\gamma W_{1},\dots ,\gamma W_{k})$ is the tuple obtained by left
translation, then $\Xi _{\ww} = \Xi _{\ww'}$.  From this it follows
that if $\gamma_{P} \in \Gamma _{P}$, then $(\gamma_{P}\gamma ) \cdot
\Xi _{\ww} = \Xi _{\ww}$.  Indeed, $\gamma_{P} \cdot \Xi _{\ww} = \Xi
_{\ww''}$, where $\ww'' = (\gamma_{P} ^{-1}W_{1},W_{2},\dots ,W_{k})$,
and any element of $\Gamma _{P}$ preserves $W_{1}$.

Now assume that $\gamma \cdot \ww $ isn't a splitting.  There are two
possibilities: (1) $\gamma W_{i}\cap \gamma W_{j}\not =\{0 \}$ for
some $i,j>1$; (2) $\gamma W_{i}\cap \gamma W_{j} =\{0 \}$ for all
$i,j>1$ and $W_{1}\cap \gamma W_{j}\not =\{0 \}$ for some $j$.  In the
first case, we have $\gamma \cdot \Xi_{\ww } = 0$ for all $\gamma $,
so the Eisenstein series is identically $0$.  In the second case, we
have $(\gamma_{P} \gamma )\cdot \Xi _{w} = 0 $ for all $\gamma_{P} \in
\Gamma _{P}$, since left translation of the tuple $\gamma_{P} \gamma
\cdot \ww $ by $\gamma_{P} ^{-1}$ preserves the incidence conditions
satisfied by the $W_{i}$.  This completes the proof.

\end{proof}

\subsection{}
For the rest of this note, we will assume that $P$ is the minimal
parabolic subgroup $P_{0}$, and will take $f\equiv 1$.  Although the
functions $E^{*}_{P,\varphi }$ are not automorphic, a certain sum of
them is.

\begin{proposition}\label{relation}
Let $W_{i}$, $i=0,\dots ,n$ be $1$-dimensional subspaces of $V$, and
let $\ww (i)$ be the tuple $(W_{0},\dots ,\hat{W}_{i},\dots ,
W_{n})$, where $\hat{W}_{i}$ means delete $W_{i}$.  Then 
\[
\varphi(\Xi _{\ww (0)}) E_{P} (f,\lambda ,g) = \sum
_{i=1}^{n} (-1)^{i+1}E_{P,\varphi }^{*} (f, \lambda , g, \ww (i)).
\]
\end{proposition}

\begin{proof}
First, the twisted series on the right are well-defined, since if
$P$ and $\ww $ are compatible then so are $P$ and $\ww (i)$ for each
$i\geq 1$.  We have the following basic relation among modular symbols
for the minimal parabolic subgroup, from \cite{ash.minimal,ash-rudolph}:
\[
\Xi _{\ww (0)} = \sum _{i=1}^{n} (-1)^{i+1} \Xi _{\ww (i)}.
\]
Note that the relations in \cite{ash.minimal,ash-rudolph} imply that
this equality holds true in $H_{d (\ww )} (Y_{\ww }, \partial Y_{\ww
}; \C )$ for \emph{any} collection of $1$-dimensional rational
subspaces $(W_{0},\dots ,W_{n})$, even with the possibility that some
$\ww (i)$ aren't splittings.  The result follows immediately from
Definition \ref{tes.def} and the fact that if $\ww' = (\gamma
W_{1},\dots ,\gamma W_{n})$ with $\gamma \in \Gamma $, then $\Xi _{\ww
'} = \Xi _{\ww (0)}$
\end{proof}

\begin{theorem}\label{main.thm}
Let $P$ be the minimal parabolic subgroup $P_{0}$, and let $\ww$ be a
compatible splitting.  Let $\varphi$ be a linear form on $H_{n-1}
(\bar Y, \partial \bar Y; \C )$.  Then the series \eqref{tes.def}
converges uniformly on compact subsets of $G\times C$, where $C$ is
the cone \eqref{coneeq}.
\end{theorem}

\begin{proof}
We begin by recalling some facts from the theory of modular symbols
associated to the minimal parabolic subgroup.  These facts are
equivalent to results in \cite{ash-rudolph}, and are just reformulated
in terms of tuples and splittings.

Let $\WWW $ be the set of all full tuples of $1$-dimensional
subspaces.  We define a function $\|\phantom{a}\|\colon \WWW
\rightarrow \Z $ as follows.  From each $1$-dimensional subspace $W$, we
choose and fix a primitive vector $v (W)\in \Z ^{n}$.  Then we set 
\[
\| \ww \| = | \det (v (W_{1}),\dots,v (W_{n}) )|.
\]

Let $\WWW _{u}\subset \WWW $ be the subset of tuples for which $\|\ww
\| = 1$.  The set $\Gamma \backslash \WWW _{u}$ is finite, where
$\Gamma $ acts by left translations.  One can show that any modular
symbol $\Xi_{\ww }$ can be written as a sum
\begin{equation}\label{modsymeq}
\Xi _{\ww } = \sum _{\ww '\in S} \Xi _{\ww '},
\end{equation}
where $S$ is a finite subset of $\WWW _{u}$ (depending on $\Xi _{\ww }$).
Moreover, the cardinality of $S$ is bounded by $p (\log \|\ww \|)$,
where $p$ is a polynomial depending only on $n$
\cite{barvinok}. 

Let $\gamma \in \Gamma $ and consider the modular symbol $\gamma \cdot
\Xi _{\ww }$.  Since $\ww $ is compatible with $P$, the space $W_{1}$
is the span of the first basis element of $V$.  Let us assume for the
moment that for $i>1$,
$W_{i}$ is the span of the
$i$th standard basis element of $V$.
This implies that $\| \gamma
\cdot \ww \|$ is the absolute value of the determinant of a fixed
$(n-1)\times (n-1)$ minor of $\gamma $.  Hence 
\begin{equation}\label{maxeq}
\| \gamma\cdot \ww \| <\!< \max \bigl\{ |\gamma _{ij}|^{n-1}\bigm | 1\leq
i,j\leq n \bigr\},
\end{equation}
where the implied constant depends only on $n$.

Let $M (\gamma )$ be the right hand side of \eqref{maxeq}.  It follows
that there is a polynomial $p_{1}$, depending only on $n$, such that 
\[
p (\log \| \gamma\cdot \ww \|) < p_{1} (\log M (\gamma )).
\]

Now consider the value $\varphi (\gamma \cdot \Xi _{\ww })$.  Since
$\Gamma \backslash \WWW _{u}$ is finite, there is a maximum value
$\varphi _{\text{max}}$ that $|\varphi| $ attains on this set.
Writing $I (\gamma , \lambda ) = \exp (\langle \rho _{P}+\lambda ,
H_{P} (\gamma g)\rangle)$, we have
\begin{equation}\label{ineq}
\sum _{\gamma \in \Gamma _{P}\backslash\Gamma } |\varphi (\gamma \cdot
\Xi _{\ww })I (\gamma , \lambda )| <\!< \sum _{\gamma \in \Gamma
_{P}\backslash\Gamma } |p_{1} (\log M (\gamma ))I (\gamma , \lambda
)|,
\end{equation}
where the implied constant depends on $n$ and $\varphi _{\text{max}}$.
The right of \eqref{ineq} has the same convergence properties as
the usual Eisenstein series, and so the proof is complete under our
assumption on $\ww $.

Now assume $W_{i}$ is a general $1$-dimensional subspace of $V$
for $i>1$.  Let $v (W_{i})_{j}$ be the $j$th coordinate of $v
(W_{i})$, and let 
\[
M (\ww ) = \max \bigl\{ |v (W_{i})_{j}|\bigm | 1\leq
i,j\leq n \bigr\}.
\]
Then  
\[
\| \gamma \cdot \ww \| <\!< M (\gamma ),
\]
where the implied constant depends on $n$ and $M (\ww) $.
The rest of the proof proceeds as above.

\end{proof}

\section{Examples}\label{examples}
\subsection{}
In this section we continue to assume that $P$ is the minimal
parabolic subgroup $P_{0}$.  We begin by discussing the connection
between the construction in this note and that of \cite{goldfeld1,
goldfeld2}.

Let $\ell$ be a positive integer, let $G = \SL _{2} (\R)$, and let $\Gamma
= \Gamma _{0} (\ell )$.  The space $X = \SL _{2} (\R)/\SO _{2} (\R)
$ is the upper halfplane $\hhh$, and we let $\hhh^{*} = \hhh \cup
\Proj^{1} (\Q )$ be the usual partial compactification obtained by adjoining
cusps.  Given a pair of cusps $(q_{1}, q_{2})$, we can determine a
full tuple $(W_{1}, W_{2})$ by setting $W_{i}$ to be the subspace of
$\Q ^{2}$ corresponding to the point $q_{i}\in \Proj^{1} (\Q )$.
Slightly abusing notation, we denote the corresponding modular symbol
by $\Xi (q_{1}, q_{2})$.

\subsection{}
We can construct an interesting linear form on the modular symbols as
follows.  Let $f$ be a fixed weight two holomorphic cuspform on
$\Gamma $.  Then we set
\[
\varphi (\Xi (q_{1},q_{2})) = -2\pi i \int _{q_{1}}^{q_{2}} f (z)\, dz,
\]
where the integration is taken along the ideal geodesic from $q_{1}$
to $q_{2}$.  Note that if $f$ is a newform, then $\varphi (\Xi (\infty
, 0))$ is the special value $-L (1,f)$.

To compute the series \eqref{tes}, let $\Gamma _{\infty } = \Gamma
\cap P$, and let $\Im \colon \hhh \rightarrow \R$ be the imaginary
part.  Let $\alpha \in \check{\aaa} _{0}$ be the standard positive
root, so that $\rho _{P} = \alpha /2$.  Write $\lambda = t\alpha $,
where $t\in \C $.  It is easy to check that $e^{\langle \lambda +\rho
_{P}, H_{P} (g)\rangle} = \Im (z)^{t+1/2}$, where $z\in \hhh $ is the
point corresponding to $g$.  Setting $(q_{1},q_{2}) = (\infty , 0)$,
we see that the corresponding tuple $\ww $ is compatible with $P$.  We
obtain
\begin{equation}\label{thiscase}
E^{*}_{P,\varphi } (\lambda ,g,\ww ) =  E^{*}_{P,\varphi } (t,z,\ww ) = \sum_{\gamma \in \Gamma _{\infty }\backslash
\Gamma } \varphi (\gamma \cdot \Xi_{\ww} ) \Im (\gamma z)^{t+1/2},
\quad t\in \C.
\end{equation}
By Theorem \ref{main.thm}, this converges for $\Re t > 1/2$.
\subsection{}\label{fivethree}
To relate this to \cite{goldfeld1, goldfeld2}, we recall the pairing
between classical modular symbols and cuspforms.  One fixes a
point $z_{0} \in \hhh^{*}$, and defines a map
\begin{align*}
[\phantom{a}]_{f}\colon &\Gamma \longrightarrow \C \\
                    &\gamma \longmapsto -2\pi i\int _{z_{0}}^{\gamma z_{0}} f
(z)\, dz.
\end{align*}
(In \cite{goldfeld1, goldfeld2}, this map is written
as $\gamma \mapsto \langle \gamma , f \rangle$.) 
One can show that this map is independent of $z_{0}$, vanishes on $\Gamma _{\infty }$, and satisfies 
\[
[\gamma \gamma ']_{f} = [\gamma ]_{f} + [\gamma ']_{f}, \quad
\text{for $\gamma , \gamma '\in \Gamma $.}
\].
Then the series in \cite{goldfeld1, goldfeld2} is defined by 
\begin{equation}\label{dgseries}
E^{*} (z, s) = \sum _{\gamma \in \Gamma _{\infty }\backslash
\Gamma } [\gamma ]_{f} \Im (\gamma z)^{s}, \quad s \in \C,
\end{equation}
which converges for $\Re s > 1$.  

To compare this with \eqref{thiscase}, let $q_{1} = \infty $ and
$q_{2} = z_{0} = 0$, and put $s = t+1/2$.  Since
\[
\int _{q_{1}}^{q_{2}} f + \int _{q_{2}}^{\gamma q_{2}} f = \int
_{q_{1}}^{\gamma q_{2}} f, 
\]
we find 
\[
\varphi (\Xi (q_{1}, q_{2})) E (z,s) + E^{*}(z,s) =
E^{*}_{P,\varphi } (s-1/2,z,\ww ),
\]
where 
\[
E (z,s) = \sum _{\gamma \in \Gamma _{\infty }\backslash \Gamma } \Im (\gamma z)^{s}
\]
is the classical nonholomorphic Eisenstein series. 

In \cite{goldfeld1, goldfeld2} it is shown that $E^{*}$ is
``automorphic up to a shift.''  Precisely, if $\gamma \in \Gamma $,
then 
\[
E^{*} (\gamma z, s) = E^{*} (z,s) - [\gamma ]_{f}E (z,s).
\]
This is easily seen to be equivalent to Proposition \ref{relation} above.

\subsection{}
Now let $G = \SL _{3} (\R)$, and let $\Gamma = \Gamma _{0} (\ell )$.  This
is the arithmetic group defined to be the subgroup of $G (\Z )$ consisting of matrices with
bottom row congruent to $(0,0,*)$ mod $\ell$.  The symmetric space $X
= \SL_{3} (\R)/\SO_{3} (\R)$
is a $5$-dimensional smooth noncompact manifold, and our modular
symbols live in $H_{2} (\bar Y, \bar \partial Y; \C )$.  

To construct an interesting linear form on these modular symbols, we
may use elements of the \emph{cuspidal cohomology} $H^{3}_{\cusp}
(\Gamma ; \C )$.  These are classes that, via the de Rham isomorphism,
correspond to $\Gamma '$-invariant differential forms $\omega = \sum
_{I} f_{I}\, d\omega _{I}$, where the coefficients are cusp forms and
$\Gamma '\subset \Gamma $ is a torsionfree subgroup of finite index.
In this context, $H^{3}_{\cusp}
(\Gamma ; \C )$ can alternatively be defined to be the kernel of the
restriction map $H^{3} (\bar Y; \C )\rightarrow H^{3} (\partial \bar
Y; \C )$.  We refer to \cite{lee.schwermer} for details.

\subsection{}
To explicitly construct classes in $H^{3}_{\cusp} (\Gamma ; \C )$ that
can be paired with modular symbols, we may use techniques of
\cite{agg}.  There it is shown that $H^{3}_{\cusp} (\Gamma ; \C )$ is
isomorphic to a space $W (\Gamma )$ of functions $f\colon \Proj
^{2}(\Z /\ell \Z )\rightarrow \C $ satisfying certain relations
\cite[Summary 3.23]{agg}.  A modular symbol $\Xi _{\ww }$ modulo
$\Gamma $ gives rise to a point $p_{\ww} \in  \Proj
^{2}(\Z /\ell \Z )$ by taking the bottom row of the matrix $(v (W_{1}),
v (W_{2}), v (W_{3}))$ \cite[Prop. 3.12]{agg}.  Hence given an element
$\alpha \in H^{3}_{\cusp} (\Gamma ; \C )$ corresponding to a function
$f_{\alpha } \in  W (\Gamma )$, we obtain a linear form by setting 
\[
\varphi (\Xi _{\ww }) = f_{\alpha } (p_{\ww }).
\] 
This linear form is induced from the intersection pairing 
\[
H_{3} (\bar Y) \times H_{2} (\partial \bar Y) \longrightarrow \C;
\]
we refer to \cite[Prop. 3.24]{agg} for details.

\subsection{}
For an explicit example, we may take $\ell = 53$.  This is the first
level for which the cuspidal cohomology is nonzero; one finds
that $\dim H^{3}_{\cusp} (\Gamma _{0} (53);\C ) = 2$.  A sample
element is given as a function in $W (\Gamma )$ in Table II
of \cite{agg}.  

To compute $E^{*}_{P,\varphi }$, we may take $\alpha \in H^{3}_{\cusp} (\Gamma _{0} (53);\C )$ to be a Hecke eigenclass.  For a
prime $p$ with $(p,53)=1$, the local $L$-factor of the
representation corresponding to $\alpha $ has the form 
\[
(1-a_{p}p^{-s }+\bar a_{p}p^{1-2s }-p^{3-3s })^{-1},
\]  
where $s\in \C $ and $a_{p}$ is the eigenvalue of a certain Hecke
operator.  If we fix an algebraic integer $\rho $ satisfying $\rho
^{2}=-11$, we find that for our Hecke eigenclass
\[
a_{2} = -2-\rho  , \quad a_{3}= -1+\rho  , \quad a_{5} = 1, \quad a_{7}=-3, \quad \dots 
\]
If we represent $\alpha $ using a function $f\in W (\Gamma )$, and
apply the formul\ae\ in \cite[Ch. V and VII]{bump}, we can obtain a
very explicit expression for $E_{P,\varphi }^{*}$.

In contrast to the $\SL_{2}$ case, the twisted Eisenstein series on
$\SL_{3}$ isn't simply automorphic up to a shift.  If we consider the
relation in Proposition \ref{relation}, we see that a certain sum of
\emph{three} twisted Eisenstein series is equal to an automorphic
function.  

\bibliographystyle{amsplain} 
\bibliography{eisenstein4}

\end{document}